\documentclass[12pt]{article}
\usepackage[ansinew]{inputenc}
\usepackage{graphicx}
\usepackage{color}
\usepackage[colorlinks]{hyperref}
\usepackage{amsmath}
\usepackage{amssymb}
\title{Primary decomposable subspaces of $k[t]$ and Right ideals of the first Weyl algebra $A_{1}(k)$ in characteristic zero}

\author{ \begin{tabular}{ll}
M. K. KOUAKOU & A. TCHOUDJEM\\
Université de Cocody & Institut Camille Jordan\\
UFR-Mathématiques & UMR 5208\\
et Informatique &  Universit\'e  Lyon 1\\
22 BP 582 Abidjan 22 & 43 bd du 11 novembre 1918\\
 & 69622 Villeurbanne cedex\\
C\^ote d'Ivoire  & France\\
  cw1kw5@yahoo.fr & tchoudjem@math.univ-lyon1.fr
\end{tabular}
}

\begin{document}
\maketitle

In the classification of right ideals $A_{1}:=k[t, \partial]$ the first Weyl algebra over a field $k$
, R. Cannings and M.P. Holland established in [3, Theorem 0.5] a bijective correspondence between primary decomposable subspaces of $R=k[t]$ and right ideals $I$ of $A_{1}:=k[t, \partial]$ the first Weyl algebra over $k$ which have non-trivial intersection with $k[t]$:

$$\Gamma:V\longmapsto \mathcal{D}(R, V )\texttt{ , }\Gamma^{-1}:I\longmapsto I\star1$$
This theorem is a very important step in this study, after Stafford's theorem [1, Lemma 4.2]. However, the theorem had been established only when the field $k$ is an algebraically closed and of characteristic zero.

In this paper we define notion of primary decomposable subspaces of $k[t]$ when $k$ is any field of characteristic zero, particulary for $\mathbb{Q}$, $\mathbb{R}$, and we show that R. Cannings and M.P. Holland's correspondence theorem holds. Thus right ideals of $A_{1}(\mathbb{Q})$, $A_{1}(\mathbb{R})$,.. are also described by this theorem.

\section { Cannings and Holland's theorem}

\subsection{Weyl algebra in characteristic zero and differential operators}

Let $k$ be a commutative field of characteristic zero and $A_{1}:=A_{1}(k)=k[t,\partial]$
where $\partial,t$ are related by $\partial t -t\partial=1$, be the first Weyl algebra over $k$.

$A_{1}$ contains the subring $R:=k[t]$ and $S:=k[\partial]$. It is well known that $A_{1}$ is an integral domain, two-sided noetherian and since the characteristic of $k$ is zero, $A_{1}$ is hereditary (see [2]). In particular, $A_{1}$ has a quotient divison ring, denoted by $Q_{1}$. For any right (resp: left ) submodule of $Q_{1}$, $M^{\ast}$ the dual as $A_{1}$-module will be identified with the set $\{u\in Q_{1}:uM\subset A_{1}\}$ (resp:$\{u\in Q_{1}:Mu\subset A_{1}\}$) when $M$ is finitely generated (see [1]).

$Q_{1}$ contains the subrings $D=k(t)[\partial]$ and $B=k(\partial)[t]$. The elements of $D$ are $k$-linear endomorphisms of $k(t)$. Precisely, if $d=a_{n}\partial^{n}+\cdot \cdot+a_{1}\partial+a_{0}$ where $a_{i}\in k(t)$ and $h\in k(t)$, then
$$d(h):= a_{n}h^{(n)}+\cdot \cdot +a_{1}h^{(1)}+a_{0}h $$
where $h^{(i)}$ denotes the $i$-th derivative of $h$ and $a_{i}h^{(i)}$ is a product in $k(t)$. One checks that:
$$(dd')(h)=d(d'(h)) \textbf{ for }d,d'\in k(t)[\partial]\;,\; h\in k(t)$$
For $V$ and $W$ two vector subspaces of $k(t)$, we set :
$$\mathcal{D}(V,W):=\{d\in k(t)[\partial]:d(V)\subset W\}$$
$\mathcal{D}(V,W)$ is called the set of differential operators from $V$ to $W$.

Notice that $\mathcal{D}(R,V)$ is an $A_{1}$ right submodule of $Q_{1}$ and $\mathcal{D}(V,R)$ is an $A_{1}$ left submodule of $Q_{1}$. If $V\subseteq R$, one notes that $\mathcal{D}(R,V)$ is a right ideal of $A_{1}$. When $V=R$, then $\mathcal{D}(R,R)=A_{1}$.

If $I$ is a right ideal of $A_{1}$, we set
$$I\star1:=\{d(1),d\in I\}$$
Clearly, $I\star1$ is a vector subspace of $k[t]$ and $I\subseteq \mathcal{D}(R,I\star1)$.

Inclusion $A_{1}\subset k(\partial)[t]$ and $A_{1}\subset k(t)[\partial]$ show that it can be defined on $A_{1}$ two notions of degree: the degree associated to "$t$" and the degree associated to "$\partial$". Naturally, those degree notions extend to $Q_{1}$.

\subsection{Stafford's theorem}

Let $I$ be  a non-zero right ideal of $A_{1}$. By J. T. Stafford in [1, Lemma 4.2], there exist $x,e \in Q_{1}$ such that:
$$(i)\; xI\subset A_{1} \texttt{ and } xI\cap k[t]\neq \{0\}\texttt{ , } (ii)\; eI\subset A_{1} \texttt{ and } eI\cap k[\partial]\neq \{0\}$$

By $(i)$ one sees that any non-zero right ideal $I$ of $A_{1}$ is isomorphic to another ideal $I'$ such that $I'\cap k[t]\neq \{0\}$, which means that $I'$ has non-trivial intersection with $k[t]$.
We denote $\mathcal{I}_{t}$ the set of right ideals $I$ of $A_{1}$ the first Weyl algebra over $k$ such that $I\cap k[t]\neq \{0\}$

 Stafford's theorem is the first step in the classification of right ideals of the first Weyl algebra $A_{1}$.

\subsection{The bijective correspondence theorem}

Let $c$ be an algebraically closed field of characteristic zero. Cannings and Holland
have defined primary decomposable subspace $V$ of $c[t]$ as finite intersections
of primary subspaces which are vector subspaces of $c[t]$ containing a power of a maximal ideal
$m$ of $c[t]$. Since $c$ is an algebraically closed field, maximal ideals of $c[t]$ are generated by one polynomial of degree one: $m=(t-\lambda)c[t]$. So, a vector subspace $V$ of $c[t]$ is primary decomposable if:

 $$V=\bigcap_{i=1}^{n}V_{i}$$
 where each $V_{i}$ contains a power of a maximal ideal $m_{i}$ of $c[t]$.

 They have established the nice well-known bijective correspondence between primary
decomposable subspaces of $c[t]$ and  $\mathcal{I}_{t}$ by:
$$\Gamma:V\longmapsto \mathcal{D}(R, V )\texttt{ , }\Gamma^{-1}:I\longmapsto I\star1$$

Since $V=\displaystyle\bigcap_{i=1}^{n}V_{i}$ and $m_{i}=<(t-\lambda_{i})^{r_{i}}>\subseteq V_{i}$,

one has $(t-\lambda_{1})^{r_{1}}\cdot\cdot\cdot(t-\lambda_{n})^{r_{n}}k[t]\subseteq V$. So, easily one sees that
$$(t-\lambda_{1})^{r_{1}}\cdot\cdot\cdot(t-\lambda_{n})^{r_{n}}k[t]\subseteq\mathcal{D}(R, V )\cap c[t]$$
 However it is not clear that $I\star1$ must be a primary decomposable subspace of $c[t]$.

Cannings and Holland's theorem use the following result, which holds even if the field is just of characteristic zero:
\bigskip

\emph{\textbf{Lemma 1}}: Let $I\in \mathcal{I}_{t}$ and $V=I\star 1$. One has:

 $I=\mathcal{D}(R,V)$ and $I^{\ast}= \mathcal{D}(V,R)$.
\bigskip

For the proof of Cannings and Holland's theorem one can see [3].
\bigskip

We note that, since $<(t-\lambda_{i})^{r_{i}}>\subseteq V_{i}$, for any $s$ in the ring $ c+(t-\lambda_{i})^{r_{i}}c[t]$, one has :
$$s\cdot V_{i}\subseteq V_{i}$$

It is this remark which will allow us to give general definition of primary decomposable subspaces of $k[t]$ for any field $k$ of characteristic zero, not necessarily algebraically closed.

\section{Primary decomposable subspaces of $k[t]$}

Here we give a general definition of primary decomposable subspaces of $k[t]$
when $k$ is any field of characteristic zero not necessarily algebraically closed and we keep the bijective correspondence of Cannings and Holland.

\bigskip

\subsection{Definitions and examples}

\bigskip
$\bullet$-\underline{ Definitions}

Let $b, h \in R = k[t]$ and $V$ a $k$-subspace of $k[t]$. We set:
$$ O(b) = \{a \in R : a'\in bR\} \;\textbf{ and } \;O(b, h) = \{a \in R : a' + ah \in bR\}$$

where $a'$ denotes the formal derivative of $a$.
$$S(V ) = \{a \in R : aV \subseteq V \}\;\textbf{ and }\; C(R, V ) = \{a \in R : aR \subseteq V \}$$
Clearly $O(b)$ and $S(V )$ are $k$-subalgebras of $k[t]$. If $b\neq 0$, the Krull dimension of $O(b)$ is $dim_{K}(O(b))=1$. The set $C(R, V )$ is an ideal of $R$
contained in both $S(V )$ and $V$ .

$\bullet$ A $k$-vector subspace $V$ of $k[t]$ is said to be primary decomposable if
$S(V )$ contains a $k$-subalgebra $O(b)$, with $b\neq0$.

 \bigskip

$\bullet$- \underline{Examples}

 $\circ$ Easily one sees that $O(b)\subseteq S(O(b,h))$ and $C(R, O(b,h))=C(R,O(b))$, in particular $O(b,h)$ is a primary decomposable subspace when $b\neq 0$.

\bigskip

Following lemmas and corollary show that classical primary decomposable subspaces are primary decomposable in the new way.

\bigskip

\emph{\textbf{Lemma 2}}: Let $k$ be a field of characteristic zero and
$\lambda_{1}$, ..,$\lambda_{n}$ finite distinct elements of $k$. Suppose that $V_{1}$ , ..,$V_{n}$  are $k$-vector subspaces $k[t]$, each $V_{i}$ contains $(t-\lambda_{i})^{r_{i}}k[t]$ for some $r_{i}\in  \mathbb{N}^{\ast}$. Then
$$O((t - \lambda_{1})^{r_{1}-1}\cdot\cdot\cdot (t - \lambda_{n})^{r_{n}-1})\subseteq S(\bigcap_{i=1}^{n}V_{i})$$

\textbf{Proof:} : One has $O((t -\lambda_{i})^{r_{i}-1}) = k + (t - \lambda_{i})^{r_{i}}k[t]$ and
$$ O((t -\lambda_{i})^{r_{i}-1}\cdot\cdot (t -\lambda_{n})^{r_{n}-1}) =\bigcap_{i=1}^{n}O((t -\lambda_{i})^{r_{i}-1})$$

An immediate consequence of this lemma is:

\bigskip

\emph{\textbf{Corollary 3}}: In the above hypothesis of lemma 2, let $$V =\bigcap_{i=1}^{n}V_{i}$$
. If $q \in C(R, V )$, then $O(q) \subseteq S(V )$.

\bigskip

\textbf{Proof:}: First one notes that if $q\in pk[t]$, then $O(q)\subseteq O(p)$. Let $b=(t-\lambda_{1})^{r_{1}}\cdot\cdot(t-\lambda_{n})^{r_{n}}$.

In the above hypothesis, one has
$$C(R,V)=\bigcap_{i=1}^{n}C(R,V_{i})=\bigcap_{i=1}^{n}(t-\lambda_{i})^{r_{i}}k[t]=(\prod_{i=1}^{n}(t-\lambda_{i})^{r_{i}})k[t]=bk[t]$$
Since $b\in (t-\lambda_{1})^{r_{1}-1}\cdot\cdot(t-\lambda_{n})^{r_{n}-1}k[t]=b_{0}k[t]$, one has $O(b_{0})V_{i}\subseteq V_{i}$ for all $i$, so
$$O(b_{0})\subseteq S(V)\texttt{ and }O(q)\subseteq O(b)\subseteq O(b_{0})$$

\bigskip

$\circ$ \underline{An opposite-example:}
\bigskip

Suppose the field $k$ is of characteristic zero and one can find $q\in k[t]$ such that: $q$ is irreducible and $deg(q)\geq 2$. Then
the vector subspace $V=k+qk[t]$ is not primary decomposable.

\subsection{Classical properties of primary decomposable subspaces }

\bigskip
Here we prove that when the field $k$ is algebraically closed of characteristic zero, those two definitions are the same.

\bigskip

\emph{\textbf{Lemma 4}}: Let $k$ be an algebraically closed field of characteristic zero and $V$ be a $k$-vector
subspace of $k[t]$ such that $S(V )$ contains a $k$-subalgebra $O(b)$ where $b\neq0$. Then $V$ is a finite intersections of subspaces which contains a power of a maximal ideal of $k[t]$.

\textbf{Proof}: Since $k$ is algebraically closed field and $b\neq0$, one can suppose

$b= (t-\lambda_{1})^{r_{1}}\cdot\cdot(t-\lambda_{n})^{r_{n}}$. Let $b^{\ast}=(t-\lambda_{1})\cdot\cdot(t-\lambda_{n})$. One has
$$ O(b) =\bigcap_{i=1}^{n}(k+(t -\lambda_{i})^{r_{i}+1}R)$$
If we suppose that $V$ is not contained in any ideal of $R$, one has $V.R=R$. Clearly
$$bb^{\ast}R=\prod_{i=1}^{n}(t -\lambda_{i})^{r_{i}+1}R\subset O(b)$$
so $bb^{\ast}R=(bb^{\ast})(RV)=(bb^{\ast}R)V=bb^{\ast}R\subset V\; (1)$. One also has
$$O(b)\cap(t-\lambda_{i})R\neq  O(b)\cap(t-\lambda_{j})R\textbf{ for all }\;i\neq j $$
in particular one has $$O(b)=[O(b)\cap(t-\lambda_{i})R]^{r_{i}+1}+ [O_{k}(b)\cap(t-\lambda_{j})R]^{r_{j}+1}\quad (2)$$
With $(1)$ and $(2)$ one gets inductively:
$$ V =\bigcap_{i=1}^{n}(V+(t -\lambda_{i})^{r_{i}+1}R)\qquad \diamondsuit$$
One also obtains usual properties of primary decomposable subspaces.

\bigskip

\emph{\textbf{Lemma 5}}: Let $k$ be a field of characteristic zero, $V$ and $W$ be  primary decomposable subspaces of $k[t]$

(1)  then $V+W$ and $V\cap W$ are  primary decomposable subspaces.

(2) If $q\in k(t)$ such that $qV\subseteq k[t]$, then $qV$ is a primary decomposable subspace.

\bigskip

\textbf{Proof}: One notes that $O(ab)\subseteq O(a)\cap O(b)$ for all $a$, $b\in k[t]$.

\bigskip

Let us recall basic properties on the subspace $O(a,h)$.

\emph{\textbf{Lemma 6 }}:

(1) $O(a)\subseteq S(O(a,h))$

(2) $C(R,O(a))=C(R,O(a,h))$

(3) $a^{2}k[t]\subset O(a)\cap O(a,h)$

(4) $\mathcal{D}(R,O(a,h))=A_{1}\cap (\partial+h)^{-1}aA_{1}$

(5) the subspace $O(a,h)$ is not contained in any proper ideal of $R$.

(6) For all $q\in O(a,h)$ such that $hcf(q,a)=1$, one has
$$O(a,h)=qO(a)+ C(R,O(a))$$

\bigskip

\textbf{Proof}: One obtains (1), (2), (3), (4) by a straightforward calculation.

Suppose $O(a,h)\subseteq gk[t]$. Then $\mathcal{D}(R,O(a,h))\subseteq gA_{1}$, and applying the $k$-automorphism $\sigma\in Aut_{k}(A_{1})$ such that $\sigma(t)=t$ and $\sigma(\partial)=\partial-h$, one obtains $\mathcal{D}(R,O(a))\subseteq gA_{1}$. Clearly the element $f=\partial^{-1} a \partial^{m+1}$ where $deg_{t}(a)=m$ belongs to $\mathcal{D}(R,O(a))=A_{1}\cap \partial^{-1}aA_{1}$. When one writes $f$ in extension, one gets exactly
$$f=a\partial^{m}+a_{m-1}\partial^{m-1}+\cdot\cdot+a_{1}\partial+(-1)^{m}m!$$
Since $f\in gA_{1}$, $(-1)^{m}m!$ must belong to $gR$. Hence $g\in k^{\star}$ and one gets (5).

Let $q$ be an element of $O(a,h)$ such that $hcf(q,a)=1$. One has also $hcf(q,a^{2})=1$, and by Bezout theorem there exist $u$, $v\in k[t]$ such that:
$$uq+va^{2}=1\quad (\ast)$$
The inclusion $qO(a)+ C(R,O(a))\subseteq O(a,h)$ is clear since $q\in O(a,h)$ and one has properties (1) and (2). Conversely let $p\in O(a,h)$. Using $(\ast)$, one gets
$$p=(pu)q+a^{2}pv\quad (\ast\ast)$$
One notes that $p(uq)=p-pva^{2}\in O(a,h)$, so $(p(uq))'+(p(uq))h \in aR$. One has $(p(uq))'+(p(uq))h=p'(uq)+p(uq)'+p(uq)h=p(uq)'+uq(p'+ph)$.

Since $q$ is chosen in $O(a,h)$, one has $p'+ph\in aR$. Then $q(up)'\in aR$, and at the end, because of $hcf(q,a)=1$, it follows that $(up)'\in aR$. Now, $up\in O(a)$ and $(\ast\ast)$ shows that $p\in qO(a)+ C(R,O(a))$.

\bigskip

\emph{\textbf{Proposition 7}}: Let $k$ be a field of characteristic zero and $V$ a $k$-vector
subspace of $k[t]$ such that $S(V )$ contains a $k$-subalgebra $O(b)$. Then
$$\mathcal{D}(R, V ) \star 1 = V$$

\textbf{Proof} :

$\bullet$ Suppose $V=O(b)$. One has $\mathcal{D}(R, O(b)) = A_{1} \cap \partial^{-1}bA_{1}$.

Suppose $b = \beta_{0} + \beta_{1}t + \cdot\cdot +\beta_{m}t^{m} \textbf{ , } \beta_{m} \neq 0$. Then $f = \partial^{-1}b\partial^{m+1} \in A_{1} \cap \partial^{-1}bA_{1}$. Let us show that $f(R) = O(b)$.
For an integer $0 \leq p \leq m$, one has:
$$\partial^{-1}t^{p}\partial^{m+1} = (t\partial - 1) \cdot(t\partial - 2) \cdot\cdot\cdot (t\partial - p)\partial^{m-p}$$
and so
$$f =\beta_{0}\partial^{m}+\sum_{p=1}^{m}  \beta_{p}(t\partial - 1) \cdot (t\partial - 2) \cdot\cdot\cdot (t\partial - p)\partial^{m-p}$$
In particular one sees that:

$(1)\; f(1) = \beta_{m}(-1)^{m}m! \neq 0$

$(2)\; f(t^{j}) = 0\;\textbf{ if }\;1 \leq j <m $

$(3)\;f(t^{m})=\beta_{0}m!$

$(4)\; deg(f(t^{j})) = j\; \textbf{when} \;j \geq m + 1$

It follows that
$$dim\dfrac{R}{f(R)}= m = dim\dfrac{R}{O(b)}$$
and since $f(R) \subseteq O(b)$, one gets $f(R) =O(b)$

$\bullet\bullet$ Suppose that $O(b)\subseteq S(V)$. One has $VO(b) = V$ and then
$$[V \mathcal{D}(R, O(b))]\star1=V [\mathcal{D}(R, O(b))\star1]=VO(b) = V$$
 By lemma 1 the equality $V\mathcal{D}(R, O(b))= \mathcal{D}(R, V)$ holds, so

$\mathcal{D}(R, V)\star1=V$.

\bigskip
Next theorem is the main result of this paper.

\bigskip

\emph{\textbf{Theorem 8}}: Let $k$ be a field of characteristic zero and $V$ a $k$-vector
subspace of $k[t]$ such that: $C(R, V ) = qk[t]$ with $q \neq 0$ and $\mathcal{D}(R, V )\star1 = V$ .
Then $S(V )$ contains some $k$-subalgebra $O(b)$ with $b\neq 0$.

\bigskip

\textbf{Proof}: One has $qk[t]\subseteq V$, and there exist $v_{0}$,$v_{1}$,...,$v_{m}$ in $V$ such that
$$V=<v_{0},v_{1},\cdot\cdot,v_{m}>\oplus\; qk[t]$$
where $<v_{0},v_{1},\cdot\cdot,v_{m}>$ denotes the vector subspace of $V$ generate by $\{v_{0},v_{1},\cdot\cdot,v_{m}\}$.
For each $v_{i}$, there exist $f_{i}\in\mathcal{D}(R, V )$ such that $f_{i}(1)=v_{i}$. Let $r=max\{deg_{\partial}(f_{i}),\;0\leq i\leq m\}$, we prove that $O(q^{r})\cdot V\subseteq V$.

Since the ideal $qk[t]$ of $R=k[t]$ is contained in $V$, we have only to prove that:
$$O(q^{r})\cdot v_{i}\subseteq V\quad \forall 0\leq i\leq m $$
We need the following lemma

\bigskip
\emph{\textbf{Lemma 9}}: Let $d=a_{p}\partial^{p}+\cdot\cdot+a_{1}\partial+a_{0}\in A_{1}(k)$ where $p\in \mathbb{N}$
, $b\in k[t]$ and $s\in O(b^{p})$ . Then $[d,s]=d\cdot s- s\cdot d\in bA_{1}$.
\bigskip

\textbf{Proof}: One has $[d,s]=[d_{1}\partial,s]=[d_{1},s]\partial+d_{1}[\partial,s]$, where $d_{1}\in A_{1}$ and $d=d_{1}\partial+a_{0}$. By induction on the $\partial$-degree of $d$, one has $[d_{1},s]\partial\in bA_{1}$. Since $deg_{\partial}(d_{1})=p-1$, it is also clear that $d_{1}b^{p}\in bA_{1}$. Finally $[d,s]\in bA_{1}$.

\bigskip

By lemma 9 above, one has $f_{i}\cdot s\in \mathcal{D}(R, V )$ and $[f_{i},s]\in qA_{1}$ for each $i$.
$$s\cdot v_{i}= s\cdot(f_{i}(1))=(s\cdot f_{i})(1)=(f_{i}\cdot s+[f_{i},s])(1)$$
One has $(f\cdot s)(1)\in V$, $[f_{i},s](1)\in qk[t]$ , it follows that $s\cdot v_{i}\in V$ and that ends the proof of theorem 8.
\bigskip

Next lemma justify the definition we gave for primary decomposable subspaces.

\bigskip

\emph{\textbf{Lemma 10}}: Let $k$ be a field of characteristic zero and suppose there exist $q$ an irreducible
element of $k[t]$ with $deg(q)\geq 2$. If $V = k + qk[t]$, then
$\mathcal{D}(R, V ) = qA_{1}$.\quad In particular $V$ is not primary decomposable subspace.

\bigskip

\textbf{Proof} : Since $q$ is irreducible, one shows by a straightforward calculation that the right ideal $qA_{1}$ is maximal. Clearly one has
$qA_{1} \subseteq \mathcal{D}(R, V )$, and $\mathcal{D}(R, V ) \neq A_{1}$ since $1 \not\in \mathcal{D}(R, V )$. So one has $qA_{1} = \mathcal{D}(R, V )$.

\bigskip

References

[1]- J. T. Stafford, ” Endomorphisms of Right Ideals of The Weyl Algebra”,
Trans. Amer. Math. Soc, 299 (1987), 623-639.

[2]- S. P. Smith and J. T. Stafford,”Differential Operators on an affine
Curves”, Proc. London. Math. Soc, (3)56 (1988), 229-259.

[3]- R. C. Cannings and M. P. Holland, ”Right Ideals of Rings of Differential
Operators”, J. Algebra, 1994, vol. 167, pp. 116-141.

\end{document}